\documentclass[11pt,reqno]{amsart}

\usepackage[T1]{fontenc}
\usepackage[utf8]{inputenc}
\usepackage{lmodern}
\usepackage{microtype}
\usepackage{amssymb,amsmath,amsfonts,amsthm,mathtools,mathrsfs}
\usepackage{aliascnt}
\usepackage[dvipsnames]{xcolor}
\usepackage[margin=1in]{geometry}
\usepackage[shortlabels]{enumitem}
\usepackage{booktabs,tabularx,array}
\usepackage{tikz}
\usepackage[colorlinks=true]{hyperref}
\hypersetup{
  linkcolor=blue!45!black,
  citecolor=magenta,
  urlcolor=blue!45!black,
  pdftitle={An entropy bound for the attractive logarithmic gas},
  pdfauthor={Matias G. Delgadino, Rishabh Gvalani, Matthew Rosenzweig}
}
\usepackage[nameinlink,capitalise,noabbrev]{cleveref}

\allowdisplaybreaks
\numberwithin{equation}{section}
\setlist{listparindent=\parindent,parsep=0pt}
\setlist[itemize]{topsep=4pt,itemsep=2pt,parsep=1pt,leftmargin=1.8em}
\setlist[enumerate]{topsep=4pt,itemsep=3pt,parsep=1pt,leftmargin=2.2em}

\theoremstyle{plain}
\newtheorem{thm}{Theorem}[section]
\newaliascnt{prop}{thm}

\aliascntresetthe{prop}
\newaliascnt{lemma}{thm}
\newtheorem{lemma}[lemma]{Lemma}
\aliascntresetthe{lemma}
\newaliascnt{cor}{thm}

\aliascntresetthe{cor}
\newaliascnt{assumption}{thm}

\aliascntresetthe{assumption}
\theoremstyle{definition}
\newaliascnt{mydef}{thm}

\aliascntresetthe{mydef}
\newaliascnt{remark}{thm}
\newtheorem{remark}[remark]{Remark}
\aliascntresetthe{remark}
\newaliascnt{conj}{thm}

\aliascntresetthe{conj}

\Crefname{thm}{Theorem}{Theorems}
\Crefname{prop}{Proposition}{Propositions}
\Crefname{lemma}{Lemma}{Lemmas}
\Crefname{cor}{Corollary}{Corollaries}
\Crefname{assumption}{Assumption}{Assumptions}
\Crefname{mydef}{Definition}{Definitions}
\Crefname{remark}{Remark}{Remarks}
\Crefname{conj}{Conjecture}{Conjectures}
\crefformat{equation}{(#2#1#3)}

\DeclarePairedDelimiter{\abs}{\lvert}{\rvert}
\DeclarePairedDelimiter{\norm}{\lVert}{\rVert}
\DeclarePairedDelimiter{\ip}{\langle}{\rangle}

\renewcommand{\d}{\mathsf d}
\newcommand{\dd}{\mathop{}\!\mathrm d}
\newcommand{\T}{\mathbb T}
\newcommand{\R}{\mathbb R}

\newcommand{\Z}{\mathbb Z}
\newcommand{\N}{\mathbb N}
\newcommand{\E}{\mathbb E}
\newcommand{\Pb}{\mathbb P}

\newcommand{\op}{\mathrm{op}}

\newcommand{\one}{\mathbf 1}
\newcommand{\be}{\beta}

\newcommand{\g}{\mathsf g}

\newcommand{\Hr}{\mathsf H}

\newcommand{\Ec}{\mathcal E}

\newcommand{\wh}{\widehat}
\newcommand{\muunif}{\mu_{\mathrm{unif}}}
\newcommand{\gm}{\mathrm{gm}}

\newcommand{\spec}{\mathrm{s}}
\newcommand{\coll}{\mathrm{c}}

\newcommand{\XN}{X_N}

\DeclareMathOperator*{\esssup}{ess\,sup}

\title[Sharp partition function bound for the attractive logarithmic gas]{Sharp partition function bound for the attractive logarithmic gas}
\author[M. G. Delgadino]{Matias G. Delgadino}
\address{Matias Delgadino, Department of Mathematics, The University of Texas at Austin, PMA 8.100, 2515 Speedway, Austin, TX 78712, USA}
\email{matias.delgadino@math.utexas.edu}
\author[R. Gvalani]{Rishabh Gvalani}
\address{Rishabh Gvalani, School of Mathematics, University of Edinburgh, James Clerk Maxwell Building, Peter Guthrie Tait Road, Edinburgh EH9 3FD, United Kingdom}
\email{rishabh.gvalani@ed.ac.uk}
\author[M. Rosenzweig]{Matthew Rosenzweig}
\address{Matthew Rosenzweig, Department of Mathematical Sciences, Carnegie Mellon University, 7127 Wean Hall, 5000 Forbes Avenue, Pittsburgh, PA 15213, USA}
\email{mrosenz2@andrew.cmu.edu}
\thanks{M.G.D. was supported by National Science Foundation grant DMS-2205937. R.G. is partially funded by the Deutsche Forschungsgemeinschaft (DFG, German Research Foundation) through SPP 2410, Hyperbolic Balance Laws in Fluid Mechanics: Complexity, Scales, Randomness (CoScaRa). M.R. was supported by National Science Foundation grants DMS-2441170, DMS-2345533, and DMS-2342349.}
\date{Working draft, November 10, 2026}
\begin{document}
\begin{abstract}
 For the attractive logarithmic potential, we prove that the partition function is bounded uniformly in the number of particles on every compact inverse-temperature interval $[0,\be_0]$ with $\be_0<B$, where $B\in(0,2\d]$ is the sharp constant in the logarithmic Hardy--Littlewood--Sobolev inequality.  The proof separates a smooth bulk, which inherits the entropy inequality, from a short-range singular interaction.  The singular contribution is controlled on the event that each box of a fixed grid contains few particles.  The complementary overcrowding event has exponentially small probability, which compensates a crude linear-in-$N$ bound on the log partition function.  The resulting estimate supplies the partition bound that can be used recycled in a few different results.  A separate Lean formalization verifies the conditional partition-bound theorem under the log-HLS assumption.
\end{abstract}
\maketitle

\section{Introduction}
\label{sec:introduction}

\subsection{The attractive logarithmic gas and the main result}
\label{subsec:intro-gas}
Fix an integer $\d\ge2$, let $\T^\d=\R^\d/\Z^\d$, and let $\muunif$ be normalized Lebesgue measure.  The mean-zero periodic logarithmic kernel $\g$, with Fourier normalization \eqref{eq:kernel-normalization}, satisfies
\begin{equation}
 \g(x)=-\log|x|+O_\d(1)
 \qquad\text{as }x\to0.
 \label{eq:intro-kernel-local-form}
\end{equation}
At inverse temperature $\be\ge0$, the attractive gas has partition function
\begin{equation}
 Z_{N,\be}
 =\int_{(\T^\d)^N}
 \exp\left\{\frac{\be}{2N}\sum_{1\le i\ne j\le N}\g(x_i-x_j)\right\}
 \dd\muunif^{\otimes N}(\XN),
 \label{eq:intro-partition-function}
\end{equation}
with $\be$ fixed as $N\to\infty$.  Our main result is a uniform bound: if every finite-entropy probability density has absolutely integrable centered interaction energy and satisfies log-HLS with a constant $B\in(0,2\d]$, then, for every $0\le\be_0<B$, one has $Z_{N,\be}\le C(\d,B,\be_0)$ for all $N\ge2$ and $0\le\be\le\be_0$.  The assumptions and conclusion are stated precisely in \Cref{thm:principal}.  The issue is to control the attractive singularity without losing uniformity in the number of particles.

In dimension two, $\g=2\pi G$ for the mean-zero Green function of $-\Delta$.  The same Gibbs weight occurs in the equilibrium formulation of periodic Patlak--Keller--Segel dynamics \cite{ChodrondeCourcelRosenzweigSerfaty2025} and in the one-species point-vortex ensemble at negative physical inverse temperature \cite{Onsager1949,Grotto2026}.

\subsection{The free-energy condition}
\label{subsec:intro-thresholds}
The macroscopic free energy is
\begin{equation}
 \Ec_\be(\mu)
 \coloneqq\Hr(\mu\mid\muunif)
 -\frac\be2\iint_{(\T^\d)^2}\g(x-y)\,\dd\mu(x)\dd\mu(y),
 \label{eq:intro-macroscopic-free-energy}
\end{equation}
where $\Hr$ denotes relative entropy and the energy domain is specified in \eqref{eq:energy-domain}.  Since $\Ec_\be(\muunif)=0$, the log-HLS assumption says that the uniform measure is a global minimizer at inverse temperature $B$.  Following \cite{ChodrondeCourcelRosenzweigSerfaty2025}, define
\begin{equation}
 \beta_{\gm}
 \coloneqq\sup\left\{b\in[0,2\d]:
 \Ec_b(\mu)\ge0\text{ for every probability measure $\mu$ with }\Hr(\mu\mid\muunif)<\infty\right\}.
 \label{eq:intro-global-minimization-threshold}
\end{equation}
Under the stated energy-domain assumption, the theorem therefore gives order-one normalization on every compact inverse-temperature interval below $\beta_{\gm}$.  It requires no identification of this threshold; the collision and spectral obstructions, and the dimension-dependent values, are discussed in \Cref{rem:threshold-context}.

The bound is the input used in our full-length paper \cite{DelgadinoGvalaniRosenzweigAttractiveI2026} for determinant asymptotics and marginal expansions.  We isolate its proof here.

\subsection{Proof idea}
\label{subsec:intro-strategy}
\begin{figure}[htbp]
\centering
\begin{tikzpicture}[x=.95cm,y=.95cm,font=\scriptsize]
\definecolor{particle}{RGB}{35,67,104}
\definecolor{overfull}{RGB}{190,50,40}
 \begin{scope}[shift={(0,0)}]
  \node at (2,4.45) {On $\mathcal A_{N,J,\delta}$: every box has $n_a\le\delta N$};
  \draw[black!40,step=1] (0,0) grid (4,4);
  \draw[black!65,line width=.5pt] (0,0) rectangle (4,4);
  \fill[particle] (0.38,0.26) circle (1.9pt);
  \fill[particle] (0.61,1.20) circle (1.9pt);
  \fill[particle] (0.19,1.51) circle (1.9pt);
  \fill[particle] (0.18,2.45) circle (1.9pt);
  \fill[particle] (0.20,3.21) circle (1.9pt);
  \fill[particle] (0.45,3.73) circle (1.9pt);
  \fill[particle] (1.24,0.31) circle (1.9pt);
  \fill[particle] (1.59,0.81) circle (1.9pt);
  \fill[particle] (1.55,1.43) circle (1.9pt);
  \fill[particle] (1.83,2.18) circle (1.9pt);
  \fill[particle] (1.25,2.23) circle (1.9pt);
  \fill[particle] (1.37,3.72) circle (1.9pt);
  \fill[particle] (2.28,0.56) circle (1.9pt);
  \fill[particle] (2.60,1.41) circle (1.9pt);
  \fill[particle] (2.19,1.29) circle (1.9pt);
  \fill[particle] (2.63,2.45) circle (1.9pt);
  \fill[particle] (2.37,2.56) circle (1.9pt);
  \fill[particle] (2.47,3.36) circle (1.9pt);
  \fill[particle] (3.71,0.64) circle (1.9pt);
  \fill[particle] (3.32,0.55) circle (1.9pt);
  \fill[particle] (3.52,1.76) circle (1.9pt);
  \fill[particle] (3.66,2.35) circle (1.9pt);
  \fill[particle] (3.84,3.23) circle (1.9pt);
  \fill[particle] (3.44,3.68) circle (1.9pt);
  \node[anchor=north east,inner sep=1.5pt,text=black!55] at (1,1) {$1$};
  \node[anchor=north east,inner sep=1.5pt,text=black!55] at (1,2) {$2$};
  \node[anchor=north east,inner sep=1.5pt,text=black!55] at (1,3) {$1$};
  \node[anchor=north east,inner sep=1.5pt,text=black!55] at (1,4) {$2$};
  \node[anchor=north east,inner sep=1.5pt,text=black!55] at (2,1) {$2$};
  \node[anchor=north east,inner sep=1.5pt,text=black!55] at (2,2) {$1$};
  \node[anchor=north east,inner sep=1.5pt,text=black!55] at (2,3) {$2$};
  \node[anchor=north east,inner sep=1.5pt,text=black!55] at (2,4) {$1$};
  \node[anchor=north east,inner sep=1.5pt,text=black!55] at (3,1) {$1$};
  \node[anchor=north east,inner sep=1.5pt,text=black!55] at (3,2) {$2$};
  \node[anchor=north east,inner sep=1.5pt,text=black!55] at (3,3) {$2$};
  \node[anchor=north east,inner sep=1.5pt,text=black!55] at (3,4) {$1$};
  \node[anchor=north east,inner sep=1.5pt,text=black!55] at (4,1) {$2$};
  \node[anchor=north east,inner sep=1.5pt,text=black!55] at (4,2) {$1$};
  \node[anchor=north east,inner sep=1.5pt,text=black!55] at (4,3) {$1$};
  \node[anchor=north east,inner sep=1.5pt,text=black!55] at (4,4) {$2$};
 \end{scope}
 \begin{scope}[shift={(5.6,0)}]
  \node at (2,4.45) {On $\mathcal B_{N,J,\delta}$: some box has $n_a>\delta N$};
  \fill[overfull!18] (1,1) rectangle (2,2);
  \draw[black!40,step=1] (0,0) grid (4,4);
  \draw[black!65,line width=.5pt] (0,0) rectangle (4,4);
  \fill[particle] (0.26,0.49) circle (1.9pt);
  \fill[particle] (0.18,1.62) circle (1.9pt);
  \fill[particle] (0.69,2.55) circle (1.9pt);
  \fill[particle] (0.76,3.37) circle (1.9pt);
  \fill[particle] (1.64,0.57) circle (1.9pt);
  \fill[particle] (1.56,1.47) circle (1.9pt);
  \fill[particle] (1.62,1.82) circle (1.9pt);
  \fill[particle] (1.19,1.64) circle (1.9pt);
  \fill[particle] (1.27,1.23) circle (1.9pt);
  \fill[particle] (1.42,1.76) circle (1.9pt);
  \fill[particle] (1.82,1.26) circle (1.9pt);
  \fill[particle] (1.53,1.17) circle (1.9pt);
  \fill[particle] (1.83,1.53) circle (1.9pt);
  \fill[particle] (1.24,1.40) circle (1.9pt);
  \fill[particle] (1.34,2.21) circle (1.9pt);
  \fill[particle] (1.22,3.50) circle (1.9pt);
  \fill[particle] (2.65,0.46) circle (1.9pt);
  \fill[particle] (2.31,1.44) circle (1.9pt);
  \fill[particle] (2.58,2.62) circle (1.9pt);
  \fill[particle] (2.67,3.74) circle (1.9pt);
  \fill[particle] (3.62,0.23) circle (1.9pt);
  \fill[particle] (3.74,1.36) circle (1.9pt);
  \fill[particle] (3.55,2.41) circle (1.9pt);
  \fill[particle] (3.67,3.29) circle (1.9pt);
  \node[anchor=north east,inner sep=1.5pt,text=black!55] at (1,1) {$1$};
  \node[anchor=north east,inner sep=1.5pt,text=black!55] at (1,2) {$1$};
  \node[anchor=north east,inner sep=1.5pt,text=black!55] at (1,3) {$1$};
  \node[anchor=north east,inner sep=1.5pt,text=black!55] at (1,4) {$1$};
  \node[anchor=north east,inner sep=1.5pt,text=black!55] at (2,1) {$1$};
  \node[anchor=north east,inner sep=1.5pt,text=overfull] at (2,2) {$9$};
  \node[anchor=north east,inner sep=1.5pt,text=black!55] at (2,3) {$1$};
  \node[anchor=north east,inner sep=1.5pt,text=black!55] at (2,4) {$1$};
  \node[anchor=north east,inner sep=1.5pt,text=black!55] at (3,1) {$1$};
  \node[anchor=north east,inner sep=1.5pt,text=black!55] at (3,2) {$1$};
  \node[anchor=north east,inner sep=1.5pt,text=black!55] at (3,3) {$1$};
  \node[anchor=north east,inner sep=1.5pt,text=black!55] at (3,4) {$1$};
  \node[anchor=north east,inner sep=1.5pt,text=black!55] at (4,1) {$1$};
  \node[anchor=north east,inner sep=1.5pt,text=black!55] at (4,2) {$1$};
  \node[anchor=north east,inner sep=1.5pt,text=black!55] at (4,3) {$1$};
  \node[anchor=north east,inner sep=1.5pt,text=black!55] at (4,4) {$1$};
 \end{scope}
 \node at (2,-.4) {$\T^\d$ with $M=J^\d$ boxes of side $2r$};
 \node at (7.6,-.4) {$\Pb(\mathcal B_{N,J,\delta})\le Me^{-N\mathsf D(\delta\Vert M^{-1})}$};
\end{tikzpicture}
\caption{Counting particles in the boxes of a fixed grid ($\d=2$, $J=4$, $N=24$, $\delta=1/8$, so $\delta N=3$).  The number in each box is its count $n_a$.  Left: on the good event $\mathcal A_{N,J,\delta}$, each box, and hence each block of adjacent boxes, carries a small fraction of the particles, so every local singular interaction has small local inverse temperature.  Right: the small set is the overcrowding event $\mathcal B_{N,J,\delta}$, on which some box exceeds $\delta N$ particles.  Since each $n_a$ is binomial with mean $N/M<\delta N$, its probability is exponentially small in $N$, which pays for the crude bound $Z_{N,\be}\le A_b^{N-1}$ used there.}
\label{fig:box-counts}
\end{figure}
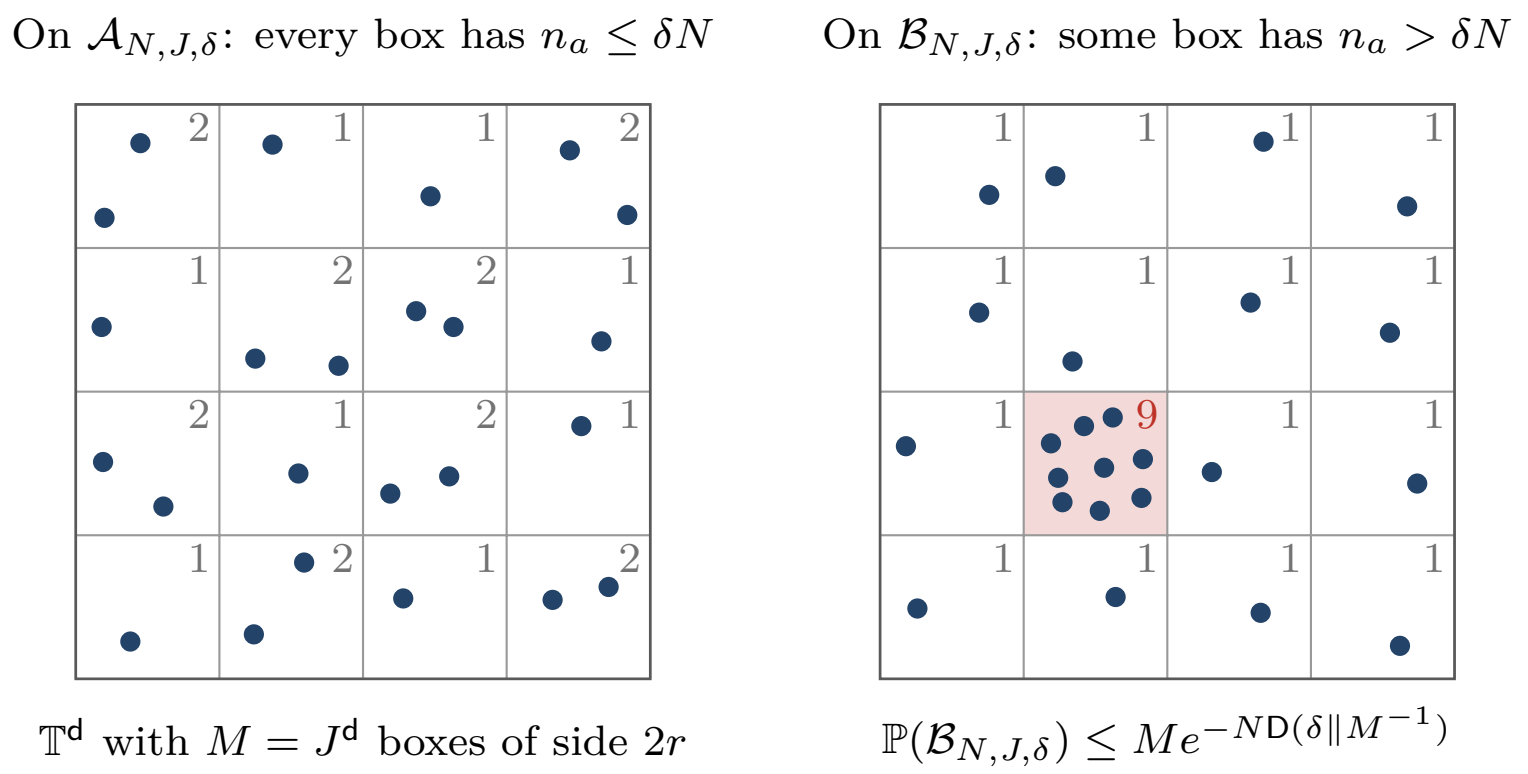

A positive-definite cutoff separates $\g$ into a smooth bulk and a centered singular part.  The bulk inherits log-HLS and is controlled by an entropy--energy estimate.  The raw singular kernel is supported near the diagonal.  After box averaging, its remaining centered interaction couples only particles in the same or adjacent boxes.  If such a block contains $m$ particles, its mean-field coefficient becomes
\begin{equation}
 \frac\be N S=\frac{\be m/N}{m}S,
 \label{eq:intro-local-temperature}
\end{equation}
so its local inverse temperature is $\be m/N$.  On the good event $\mathcal A_{N,J,\delta}$ that every box contains at most $\delta N$ particles, a perturbative estimate for centered arrays controls all local blocks; see Figure~\ref{fig:box-counts}.  The small set is its complement $\mathcal B_{N,J,\delta}$, the overcrowding event that some box contains more than $\delta N$ particles.  Since the box counts are binomial with mean $N/M<\delta N$, this event has exponentially small probability, which compensates a crude exponential bound on the Gibbs weight.  The grid is fixed independently of $N$.

\subsection{Relation to previous work}
\label{subsec:intro-related-work}
Chodron de Courcel, Rosenzweig, and Serfaty \cite[Section~8.2]{ChodrondeCourcelRosenzweigSerfaty2025} estimate the partition function by coarse graining the empirical measure, using the modulated-free-energy framework of Bresch, Jabin, and Wang \cite{BreschJabinWang2023}.  Here, the cutoff preserves the macroscopic inequality for the bulk while local particle counts control the singular remainder.  The bulk estimate is a specialization of Ren and Wang \cite[Corollary~4 and Remark~3]{RenWang2025}; we include a cell-projection proof.  For the singular part, we use centering in the sense of Hoeffding \cite{Hoeffding1948} and a random index split in the spirit of de la Pe\~na and Montgomery-Smith \cite{delaPenaMontgomerySmith1995}.

Our companion work \cite{DelgadinoGvalaniRosenzweig2026} gives Hilbert--Schmidt partition estimates, including attractive logarithmic estimates under perturbative hypotheses.  Grotto \cite[Theorem~1]{Grotto2026} proves qualitative determinant convergence on the two-dimensional torus at sufficiently small attractive inverse temperature.  Wang and Zhao \cite[Section~1.4]{WangZhao2026} announce corresponding full-subcritical results for the planar gas at the explicit log-HLS equilibrium \cite{WangZhaoAttractive2026}.  The present theorem concerns the uniform torus reference and uses the assumed macroscopic inequality throughout its strict subcritical range.

\paragraph{Organization.}
Section~\ref{sec:setting} states the theorem and gives a proof roadmap.  Sections~\ref{sec:physical-cutoff-proof}--\ref{sec:local-singular} establish the partition bound.  Section~\ref{sec:canonical-arrays} and Appendix~\ref{app:continuous-proof} supply the centered-array and continuous-kernel estimates.

\subsection{Statement on AI use}
ChatGPT, Codex, and Claude assisted with mathematical exploration, literature searches, and exposition.  This working draft remains subject to author verification.  The separate Lean formalization concerns the conditional partition-bound theorem, as explained after \Cref{thm:principal}; accompanying AI reviews are not independent human reviews.

\section{Setting and partition bound}
\label{sec:setting}
Let $\d\ge2$ and let $\muunif$ be normalized Lebesgue measure on $\T^\d=\R^\d/\Z^\d$.  The mean-zero periodic logarithmic kernel is normalized by
\begin{equation}
 \wh\g(k)=\frac{c_\d}{(2\pi|k|)^\d}\quad(k\ne0),
 \qquad c_\d=2^{\d-1}\pi^{\d/2}\Gamma(\d/2),
 \qquad \g(x)=-\log|x|+O_\d(1)
 \quad(x\to0).
 \label{eq:kernel-normalization}
\end{equation}
Here, $\wh f(k)=\int f(x)e^{-2\pi ik\cdot x}\,\dd\muunif(x)$ and $\wh\g(0)=0$.
For independent $X_1,\ldots,X_N$ with law $\muunif$, set
\begin{equation}
 \mathsf W_N(K)\coloneqq\frac1N\sum_{i\ne j}K(X_i,X_j),
 \qquad
 Z_{N,\be}\coloneqq\E\big[e^{\frac\be2\mathsf W_N(\g)}\big],
 \label{eq:homogeneous-partition-function}
\end{equation}
where a function of one variable denotes the translation-invariant kernel $K(x-y)$.

We make the energy domain explicit: for every probability density $\rho$ of finite entropy, assume
\begin{equation}
 \iint_{(\T^\d)^2}\abs{\g(x-y)(\rho(x)-1)(\rho(y)-1)}
 \,\dd\muunif(x)\dd\muunif(y)<\infty.
 \label{eq:energy-domain}
\end{equation}
Assume also the log-HLS inequality: for some $B\in(0,2\d]$,
\begin{equation}
 \Hr(\rho\muunif\mid\muunif)
 \coloneqq\int\rho\log\rho\,\dd\muunif
 \ge\frac B2\iint\g(x-y)(\rho(x)-1)(\rho(y)-1)
                  \,\dd\muunif(x)\dd\muunif(y)
 \label{eq:strict-subcritical-macroscopic-inequality}
\end{equation}
for every probability density of finite entropy.  Both \eqref{eq:energy-domain} and \eqref{eq:strict-subcritical-macroscopic-inequality} are hypotheses of the result below; their validity is not proved here.  No value or optimality of $B$ is needed.

\begin{thm}[Uniform partition bound]
\label{thm:principal}
Under \eqref{eq:energy-domain} and \eqref{eq:strict-subcritical-macroscopic-inequality}, for every $0\le\be_0<B$ there is a finite constant $C=C(\d,B,\be_0)$ such that
\begin{equation}
 \sup_{N\ge2}\sup_{0\le\be\le\be_0}Z_{N,\be}\le C.
 \label{eq:principal-order-one}
\end{equation}
\end{thm}

\paragraph{Formalization.}
The Lean theorem \texttt{LogGasPaper.main\_with\_domain} formalizes
\Cref{thm:principal}, using the Fourier characterization of
$\g$ and the explicit hypotheses \eqref{eq:energy-domain} and
\eqref{eq:strict-subcritical-macroscopic-inequality}.
The Lean sources and checking instructions are available in the
\href{https://github.com/ZaytaM/DGR-Collab-Public/tree/loggas-partition-bound-20261110/projects/attractive-log-equilibrium}{public verification repository};
the \href{https://github.com/ZaytaM/DGR-Collab-Public/releases/tag/loggas-partition-bound-20261110}{full verification evidence}
is distributed with the release.  This formalization does not establish those hypotheses
for any positive $B$, identify the sharp threshold, or separately certify the
local asymptotic formula for $\g$.  Its scope is the conditional partition bound,
not the full manuscript or the determinant and marginal expansions of the
companion paper.

\begin{remark}[Thresholds and companion results]
\label{rem:threshold-context}
The collision threshold is $\beta_{\coll}=2\d$: for $N=2$, the Gibbs weight is $e^{b\g(x_1-x_2)/2}$ and is integrable exactly when $b<2\d$.  The spectral threshold is $\beta_{\spec}=\norm{T}_{\op}^{-1}$, where $T$ is the integral operator with kernel $\g$ on $L^2(\muunif)$; below it, the free-energy Hessian at the uniform measure is positive.  Concentration and the lowest nonzero Fourier modes give $\beta_{\gm}\le\min\{\beta_{\coll},\beta_{\spec}\}$.  Global minimality is a stronger requirement than local stability.

The value of $\beta_{\gm}$ is a question about the sharp periodic logarithmic Hardy--Littlewood--Sobolev, equivalently Beckner--Onofri, inequality.  The sharp Euclidean and spherical inequalities were established by Carlen and Loss \cite{CarlenLoss1992} and Beckner \cite{Beckner1993}, and Okikiolu \cite{Okikiolu2008} developed the compact-manifold framework relating them to regularized spectral traces.  The periodic problem has the same local singularity, but its discrete spectrum and the absence of dilation symmetry create independent collision and spectral obstructions.  In dimension two, the result of Lin and Lucia \cite[Theorem~1.1 and (1.4)]{LinLucia2007} gives $\beta_{\gm}(2)=4=\beta_{\coll}$, so \Cref{thm:principal} covers the full collision-subcritical range.  In higher dimensions, forthcoming work of Lei and the third author \cite{LeiRosenzweig2026} shows that $\beta_{\gm}(\d)=2\d$ for $2\le\d\le10$, that $\beta_{\gm}(\d)=\beta_{\spec}(\d)$ for $\d\ge12$, and that $\beta_{\gm}(11)<\beta_{\spec}(11)<22$.  \Cref{thm:principal} does not use these values; it uses the energy-integrability assumption and the inequality with some constant $B\in(0,2\d]$.

In the full-length paper \cite{DelgadinoGvalaniRosenzweigAttractiveI2026}, a product trial measure gives a strictly positive lower limit for $N^{-1}\log Z_{N,\be}$ when $\beta_{\gm}<\be<2\d$.  The same paper uses the uniform partition bound at a slightly larger inverse temperature to obtain the Carleman--Fredholm determinant limit with error of order $N^{-1}$, extensions to bi-Lipschitz images of the uniform measure, and all-order expansions of fixed marginals.  These companion results are outside the scope of the present theorem and its Lean formalization.
\end{remark}

The proof has three ingredients.  A positive-definite cutoff (Section~\ref{sec:physical-cutoff-proof}) writes $\g=b_r+\widetilde s_r$, with a smooth bulk $b_r$ and a centered singular part $\widetilde s_r=s_r-c_r$, both positive semidefinite.  The raw kernel $s_r$ is supported within distance $r$ of the diagonal; subtracting its mean $c_r$ removes this support property.  The bulk inherits \eqref{eq:strict-subcritical-macroscopic-inequality} and is controlled by a continuous-kernel entropy bound (Lemma~\ref{lem:continuous-entropy-partition}).  For the singular part, box averaging separates a perturbative averaging term from a kernel centered within each box, which vanishes between nonadjacent boxes.  If no box of side $2r$ is overcrowded, each interacting block contains only a small fraction of the particles, and a moment estimate for centered arrays (Lemma~\ref{lem:canonical-array-exponential}) applies.  The overcrowded event has exponentially small probability, which compensates an exponentially growing bound for the partition function (Section~\ref{subsec:proof-bad-set}).

\textbf{Boxes and the good event.}
For an integer $J\ge3$, partition the torus into
\begin{equation}
 A_a=J^{-1}(a+[0,1)^\d),\qquad
 a\in(\Z/J\Z)^\d,\qquad M=J^\d.
 \label{eq:common-box-partition}
\end{equation}
Write $L_i=a$ when $X_i\in A_a$, and $n_a=\#\{i:L_i=a\}$.
Let $\nu_a=M\one_{A_a}\muunif$ be the uniform probability measure on $A_a$, and let $Q_J$ denote box averaging: $(Q_Ju)|_{A_a}=\int u\,\dd\nu_a$.
Both parts of the proof use this grid, at a resolution independent of $N$: the singular part at the resolution tied to the cutoff below, and the bulk part (Appendix~\ref{app:continuous-proof}) at a finer resolution chosen after $r$ has been fixed.

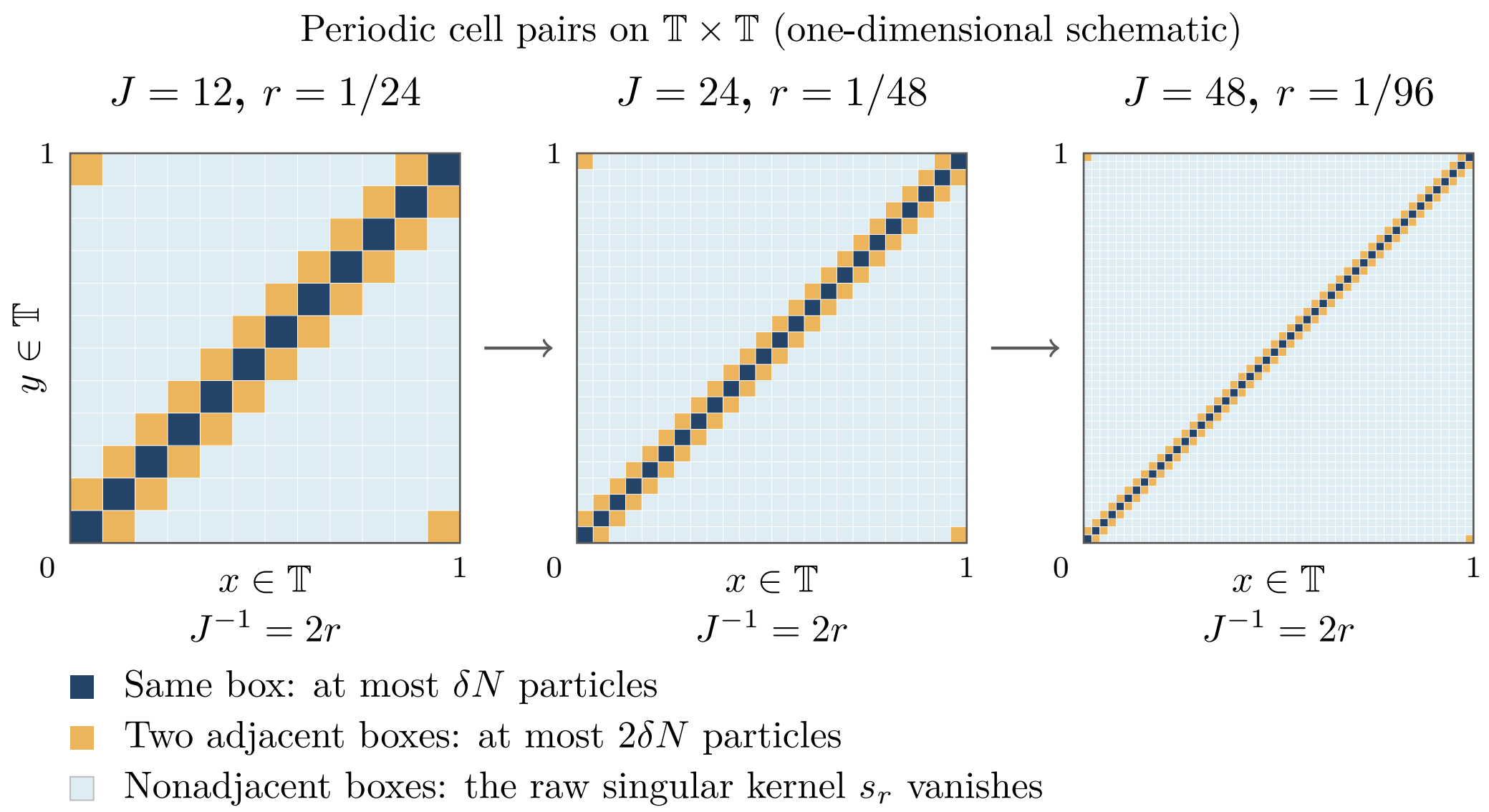
\begin{figure}[htbp]
\centering
\begin{tikzpicture}[x=.9cm,y=.9cm,font=\small]
\node at (7.2,5.25) {Periodic cell pairs on $\mathbb T\times\mathbb T$ (one-dimensional schematic)};
\definecolor{cellself}{RGB}{35,67,104}
\definecolor{cellnear}{RGB}{237,181,91}
\definecolor{cellfar}{RGB}{221,237,243}
\foreach \cells/\shift/\cutoff in {12/0/{1/24},24/5.2/{1/48},48/10.4/{1/96}} {
 \begin{scope}[shift={(\shift,0)}]
  \node[font=\bfseries] at (2,4.6) {$J=\cells$, $r=\cutoff$};
  \fill[cellfar] (0,0) rectangle (4,4);
  \pgfmathtruncatemacro{\lastcell}{\cells-1}
  \foreach \a in {0,...,\lastcell} {
   \foreach \offset in {-1,1} {
    \pgfmathtruncatemacro{\b}{mod(\a+\offset+\cells,\cells)}
    \fill[cellnear] ({4*\a/\cells},{4*\b/\cells}) rectangle ({4*(\a+1)/\cells},{4*(\b+1)/\cells});
   }
   \fill[cellself] ({4*\a/\cells},{4*\a/\cells}) rectangle ({4*(\a+1)/\cells},{4*(\a+1)/\cells});
  }
  \foreach \a in {0,...,\cells} {
   \draw[white,line width=.13pt] ({4*\a/\cells},0)--({4*\a/\cells},4);
   \draw[white,line width=.13pt] (0,{4*\a/\cells})--(4,{4*\a/\cells});
  }
  \draw[black!65,line width=.5pt] (0,0) rectangle (4,4);
  \node[below left,font=\scriptsize] at (0,0) {$0$};
  \node[below,font=\scriptsize] at (4,0) {$1$};
  \node[left,font=\scriptsize] at (0,4) {$1$};
  \node[below] at (2,-.08) {$x\in\mathbb T$};
  \ifnum\cells=12\node[rotate=90] at (-.42,2) {$y\in\mathbb T$};\fi
  \node at (2,-.85) {$J^{-1}=2r$};
 \end{scope}
}
\foreach \shift in {0,5.2} {\draw[->,thick,black!65] ({4.25+\shift},2)--({4.95+\shift},2);}
\fill[cellself] (0,-1.6) rectangle (.24,-1.36);
\node[anchor=west] at (.4,-1.48) {Same box: at most $\delta N$ particles};
\fill[cellnear] (0,-2.12) rectangle (.24,-1.88);
\node[anchor=west] at (.4,-2) {Two adjacent boxes: at most $2\delta N$ particles};
\fill[cellfar,draw=black!25] (0,-2.64) rectangle (.24,-2.4);
\node[anchor=west] at (.4,-2.52) {Nonadjacent boxes: the raw singular kernel $s_r$ vanishes};
\end{tikzpicture}

\caption{The grid is chosen to match the singular cutoff: the cell side is $2r$.  In this one-dimensional schematic, $s_r(x-y)$ can be nonzero only in a diagonal square or in one of its two neighboring squares, including the periodic corner squares.  The colored boxes cover its support; they are not all filled by it.  Refining the grid and reducing $r$ together narrows the band while keeping three possible partner boxes per box, the box itself included ($3^\d$ in dimension $\d$).}
\label{fig:periodic-cell-pair-band}
\end{figure}

For the singular part, fix $0<\delta<1$ and choose the cutoff and the occupancy event together:
\begin{equation}
 r=\frac1{2J},\qquad
 \mathcal A_{N,J,\delta}=\{n_a\le\delta N\ \text{for every }a\},
 \qquad
 \mathcal B_{N,J,\delta}=\mathcal A_{N,J,\delta}^{\mathsf c}.
 \label{eq:physical-neighbor-events}
\end{equation}
The raw singular kernel $s_r$ constructed in Section~\ref{sec:physical-cutoff-proof} is supported at distance at most $r$.
Nonadjacent boxes are separated by at least $J^{-1}=2r$, so only particles in the same box or in adjacent boxes interact through it.  Here, adjacency includes shared faces, edges, and corners, with periodic labels, so each box has $3^\d-1$ neighbors.
On $\mathcal A_{N,J,\delta}$, one box contains at most $\delta N$ particles and two boxes at most $2\delta N$; see Figure~\ref{fig:periodic-cell-pair-band}.

\section{Separate the bulk and the diagonal}
\label{sec:physical-cutoff-proof}
With $r=(2J)^{-1}$ as above, choose a nonzero, nonnegative, radial $\chi\in C_c^\infty(B_{1/2})$ and set $\phi=\chi*\chi$.  Define
\begin{equation}
 \begin{aligned}
 s_r(x)&=\frac1{\phi(0)}\sum_{m\in\Z^\d}\int_0^r
                  \phi\!\left(\frac{x+m}{t}\right)\frac{\dd t}{t},\\
 c_r&=\int s_r\,\dd\muunif,\qquad
 \widetilde s_r=s_r-c_r,\qquad b_r=\g-\widetilde s_r.
 \end{aligned}
 \label{eq:compact-positive-log-split}
\end{equation}
The raw kernel $s_r$ is nonnegative, supported in $\{\operatorname{dist}_{\T^\d}(x,0)\le r\}$, and $c_r=O(r^\d)$.  Since $\wh\phi=|\wh\chi|^2\ge0$, radial Fourier inversion gives, for $k\ne0$,
\begin{equation}
 0\le\wh{\widetilde s_r}(k)
 =\frac1{\phi(0)}\int_0^r t^{\d-1}\wh\phi(tk)\,\dd t
 \le\frac1{\phi(0)}\int_0^\infty t^{\d-1}\wh\phi(tk)\,\dd t
 =\wh\g(k).
 \label{eq:compact-split-fourier-order}
\end{equation}
Thus both centered parts are positive semidefinite.  The Fourier coefficients of $b_r$, given by the integral from $r$ to infinity, decay rapidly, so $b_r$ is smooth; it is also even and has mean zero.  Its energy form is bounded by that of $\g$, so it satisfies \eqref{eq:strict-subcritical-macroscopic-inequality} with the same $B$.

Group the label vectors $L$ by their occupancies $n=(n_a)$.  For a fixed $n$, retag the particles as $Y_{a,1},\ldots,Y_{a,n_a}$ in box $A_a$, and let $\E_n$ denote expectation under their independent product law, with each $Y_{a,p}$ distributed according to $\nu_a$.  Since the energy is symmetric in the particles,
\begin{equation}
 Z_{N,\be}
 =\sum_{\substack{n\in\N_0^M\\\sum_a n_a=N}}
   \pi_N(n)\,\E_n\big[e^{\frac\be2\mathsf W_N(\g)}\big],
 \qquad
 \pi_N(n)=\frac{N!}{M^N\prod_a n_a!}.
 \label{eq:partition-occupancy-sum}
\end{equation}
Here, $\mathsf W_N$ is evaluated on the retagged particles in any order, and the multinomial weights $\pi_N(n)$ sum to one.  Split this sum into good vectors, with $\max_a n_a\le\delta N$, and the remaining bad vectors.  Fix $0<\be_0<B$ and conjugate $p,q>1$ with $p\be_0<B$.  H\"older on the good contribution gives
\begin{equation}
 \begin{aligned}
 Z_{N,\be}
 &\le
 \Big(\E\big[e^{\frac{p\be}{2}\mathsf W_N(b_r)}\big]\Big)^{1/p}
 \Big(\E\big[\one_{\mathcal A_{N,J,\delta}}
        e^{\frac{q\be}{2}\mathsf W_N(\widetilde s_r)}\big]\Big)^{1/q}\\
 &\qquad+\E\big[\one_{\mathcal B_{N,J,\delta}}
                     e^{\frac\be2\mathsf W_N(\g)}\big],
 \qquad 0\le\be\le\be_0.
 \end{aligned}
 \label{eq:physical-cutoff-holder}
\end{equation}
The indicator stays on the singular factor, whose coupling $q\be$ can be large.  The next three sections bound the bad event, the bulk, and the singular factor.

\section{Discard overcrowded boxes}
\label{subsec:proof-bad-set}
We first need a rough bound on the full Gibbs weight.  For $0<b<2\d$, let
$A_b=\int e^{b\g/2}\,\dd\muunif<\infty$.
Convexity of the exponential, followed by integration conditional on $X_i$, gives
\begin{equation}
 \begin{aligned}
 Z_{N,b}
 &=\E\left[\exp\left\{\frac1N\sum_{i=1}^N
                   \frac b2\sum_{j\ne i}\g(X_i-X_j)\right\}\right]\\
 &\le\frac1N\sum_{i=1}^N
       \E\left[e^{\frac b2\sum_{j\ne i}\g(X_i-X_j)}\right]
 =A_b^{N-1}.
 \end{aligned}
 \label{eq:rough-partition-bound}
\end{equation}
This uses only the local logarithmic singularity.  Moreover, $A_b\ge1$ by Jensen, and H\"older gives $Z_{N,\gamma}\le Z_{N,b}^{\gamma/b}\le A_b^{N-1}$ for $0\le\gamma\le b$.

Each $n_a$ is binomial with parameters $N,M^{-1}$.  For $M^{-1}<\delta$, Chernoff's inequality and a union bound give
\begin{equation}
 \begin{aligned}
 \Pb(\mathcal B_{N,J,\delta})
 &\le M e^{-N\mathsf D(\delta\Vert M^{-1})},\\
 \mathsf D(\delta\Vert M^{-1})
 &=\delta\log(\delta M)
   +(1-\delta)\log\frac{1-\delta}{1-M^{-1}}
 \longrightarrow\infty\quad(J\to\infty).
 \end{aligned}
 \label{eq:physical-bad-probability}
\end{equation}
Choose $s>1$ with $b=s\be_0<2\d$ and $s'=s/(s-1)$.  H\"older and \eqref{eq:rough-partition-bound} yield
\begin{equation}
 \E\big[\one_{\mathcal B_{N,J,\delta}}e^{\frac\be2\mathsf W_N(\g)}\big]
 \le M^{1/s'}\exp\left\{N
 \left(\frac{\log A_b}{s}
       -\frac{\mathsf D(\delta\Vert M^{-1})}{s'}\right)\right\}.
 \label{eq:physical-bad-gibbs-bound}
\end{equation}
For fixed $\delta>0$, choose $J$ large enough to make the exponent negative, then keep this grid fixed as $N$ varies.

\section{The continuous bulk}
\label{sec:continuous-localized-reduction}
The next lemma is a continuous-kernel entropy--energy partition bound.  It is a special case of the bounded-interaction torus result of Ren and Wang \cite[Corollary~4 and Remark~3]{RenWang2025}: in their notation, take $W_+=0$, $W_-=\be K$, reference measure $\muunif$, and $\lambda=0$, with $0<\varepsilon\le1/2$ and $(1+\varepsilon)\be<B$; the diagonal terms included in their energy contribute at most $\be\|K\|_\infty/2$ to the exponent.  For completeness, Appendix~\ref{app:continuous-proof} gives a self-contained proof using cell averages, an energy projection, and the estimate of Section~\ref{sec:canonical-arrays}.

\begin{lemma}[Continuous kernels controlled by entropy]
\label{lem:continuous-entropy-partition}
Let $K\in C(\T^\d\times\T^\d)$ be real and symmetric, with $\int K(x,y)\,\dd\muunif(y)=0$.  Suppose its integral operator $T_K$ is nonnegative and, for some $B>0$,
\begin{equation}
 \Hr(\rho\muunif\mid\muunif)
 \ge\frac B2\ip{\rho-1,T_K(\rho-1)}
 \quad\text{for every finite-entropy probability density }\rho.
 \label{eq:continuous-entropy-assumption}
\end{equation}
Then, for every $0\le\be_0<B$,
\begin{equation}
 \sup_{N\ge1}\sup_{0\le\be\le\be_0}
 \E\big[e^{\frac\be2\mathsf W_N(K)}\big]<\infty.
 \label{eq:continuous-partition-conclusion}
\end{equation}
\end{lemma}

For $K(x,y)=b_r(x-y)$, the hypotheses were verified in Section~\ref{sec:physical-cutoff-proof}: $b_r$ is smooth, even, and mean zero, its operator is nonnegative by \eqref{eq:compact-split-fourier-order}, and it satisfies \eqref{eq:continuous-entropy-assumption} with the $B$ of \eqref{eq:strict-subcritical-macroscopic-inequality}.  Since $p\be_0<B$, the lemma gives
\begin{equation}
 \sup_{N\ge1}\sup_{0\le\be\le\be_0}
 \E\big[e^{\frac{p\be}{2}\mathsf W_N(b_r)}\big]<\infty
 \qquad\text{for every fixed }0<r<1/4.
 \label{eq:physical-continuous-bulk-bound}
\end{equation}

\section{The singular part and the conclusion}
\label{sec:local-singular}
The good event makes the singular factor perturbative through the number of particles in each interacting pair of boxes.  The moment estimate used below is proved in Lemma~\ref{lem:canonical-array-exponential}.
\begin{lemma}[Short-range estimate under a local count bound]
\label{lem:physical-local-count-bound}
For every $0<t_0<\infty$, there are $0<\delta_0<1$ and an integer $J_0\ge3$ such that, for each fixed $0<\delta\le\delta_0$ and integer $J\ge J_0$, with $r=(2J)^{-1}$,
\begin{equation}
 \sup_{N\ge1}\sup_{0\le t\le t_0}
 \E\Big[\one_{\mathcal A_{N,J,\delta}}
 e^{\frac t2\mathsf W_N(\widetilde s_r)}\Big]<\infty.
 \label{eq:physical-local-count-conclusion}
\end{equation}
The constant may depend on $J,\delta,t_0$; these parameters are independent of $N$.
\end{lemma}

\begin{proof}
Fix $0\le t\le t_0$.  We separate off the box averages, then use \eqref{eq:partition-occupancy-sum} to estimate the remaining integral for each good occupancy vector.

\textbf{1. Split the kernel and apply H\"older.}
Apply the box average $Q_J$ in each variable and write
\begin{equation}
 \widetilde s_r=a_r+h_r,
 \qquad h_r=(I-Q_{J,x})(I-Q_{J,y})s_r,
 \qquad a_r=Q_{J,x}s_r+Q_{J,y}s_r-Q_{J,x}Q_{J,y}s_r-c_r.
 \label{eq:physical-local-centering}
\end{equation}
By construction, $h_r$ has zero average in each variable over every box; this is the centering needed in Step~3.  Both $a_r$ and $h_r$ also have zero average in each variable under $\muunif$, which is the centering needed in Step~2.  The particle sum splits exactly as
\begin{equation}
 \frac t2\mathsf W_N(\widetilde s_r)
 =\frac tN\sum_{i<j}a_r(X_i,X_j)
  +\frac tN\sum_{i<j}h_r(X_i,X_j).
 \label{eq:physical-particle-sum-splitting}
\end{equation}
Cauchy--Schwarz therefore gives
\begin{equation}
 \E\big[\one_{\mathcal A_{N,J,\delta}}e^{\frac t2\mathsf W_N(\widetilde s_r)}\big]
 \le
 \Big(\E\big[e^{t\mathsf W_N(a_r)}\big]\Big)^{1/2}
 \Big(\E\big[\one_{\mathcal A_{N,J,\delta}}e^{t\mathsf W_N(h_r)}\big]\Big)^{1/2}.
 \label{eq:physical-local-holder}
\end{equation}
We bound the first factor by taking $r$ small and the second by taking $\delta$ small. 

\textbf{2. The averaging term is perturbative.}
Every box average of the nonnegative kernel $s_r$ is at most $Mc_r$.  Since $c_r\le Cr^\d$ and $Mr^\d=2^{-\d}$, we have
\begin{equation}
 \|a_r\|_\infty\le C,
 \qquad
 \sup_x\int|a_r(x,y)|\,\dd\muunif(y)\le4c_r\le Cr^\d.
 \label{eq:physical-average-smallness}
\end{equation}
The second estimate follows because each averaging term has integral $c_r$ in either variable.  These bounds imply
$\sup_x\int(e^{|a_r(x,y)|/\kappa}-1)\,\dd\muunif(y)
\le Cr^\d(e^{C/\kappa}-1)$.
Taking $\kappa(r)=C'/\log(1/r)$ with $C'$ large enough makes the right side at most one for small $r$.  Expanding the exponential and using symmetry gives, for every integer $p\ge1$,
\begin{equation}
 \sup_x\int|a_r(x,y)|^p\,\dd\muunif(y)
 +\sup_y\int|a_r(x,y)|^p\,\dd\muunif(x)
 \le2\kappa(r)^p p!,
 \qquad \kappa(r)\longrightarrow0.
 \label{eq:physical-average-factorial}
\end{equation}
Apply Lemma~\ref{lem:canonical-array-exponential} with $m=N$, $\kappa=\kappa(r)$, and $\theta=2t$:
\begin{equation}
 \E\big[e^{t\mathsf W_N(a_r)}\big]
 =\E\Big[e^{\frac{2t}{N}\sum_{i<j}a_r(X_i,X_j)}\Big]
 \le\frac{2}{1-(2tC_{\mathrm u}\kappa(r))^2}
 \le\frac83,
 \qquad 0\le t\le t_0,
 \label{eq:physical-average-exponential-bound}
\end{equation}
provided $J\ge J_0$, with $J_0$ chosen so that $2t_0C_{\mathrm u}\kappa((2J_0)^{-1})\le1/2$.

\textbf{3. Few particles make the singular term perturbative.}
Under the box-uniform measures $\nu_a$, there is a constant $\kappa_*$, independent of $J$, such that, for every integer $p\ge1$,
\begin{equation}
 \sup_{y\in A_b}\int_{A_a}|h_r(x,y)|^p\,\dd\nu_a(x)
 +\sup_{x\in A_a}\int_{A_b}|h_r(x,y)|^p\,\dd\nu_b(y)
 \le2\kappa_*^p p!,
 \qquad p\ge1.
 \label{eq:physical-singular-factorial}
\end{equation}
Indeed, Cauchy--Schwarz gives $0\le\phi\le\phi(0)$ in \eqref{eq:compact-positive-log-split}, so $s_r(z)\le(\log(r/\operatorname{dist}_{\T^\d}(z,0)))_+$.  Thus
$M\int s_r^p\,\dd\muunif\le M|B_r|p!/\d^p\le C^p p!$.
The averages subtracted to form $h_r$ are bounded by $Mc_r\le C$, which proves the same estimate for $h_r$ with a larger constant.

The occupancy decomposition gives exactly
\begin{equation}
 \E\big[\one_{\mathcal A_{N,J,\delta}}e^{t\mathsf W_N(h_r)}\big]
 =\sum_{\substack{n\in\N_0^M,\ \sum_a n_a=N\\\max_a n_a\le\delta N}}
   \pi_N(n)\,\E_n\big[e^{t\mathsf W_N(h_r)}\big].
 \label{eq:physical-good-occupancy-sum}
\end{equation}
Fix one vector in this sum.  Under $\E_n$, its particles form $M$ independent groups, one per box, each of size at most $\delta N$, and the particles in box $A_a$ have law $\nu_a$.  The kernel $h_r$ is centered under $\nu_a$ in each variable, which is the hypothesis \eqref{eq:canonical-array-centering}, and its moments satisfy \eqref{eq:canonical-array-kernel-moments} with $\kappa=\kappa_*$ by \eqref{eq:physical-singular-factorial}.  Only particles in the same or in adjacent boxes interact: $s_r$ vanishes on other box products, and box averaging preserves this property for $h_r$.  Thus each group interacts with at most $K-1$ others, where $K=3^\d$.

Apply the grouped conclusion \eqref{eq:canonical-grouped-exponential} of Lemma~\ref{lem:canonical-array-exponential} with $m=N$, $m_*=\max_a n_a$, and $u=2t$, since $\mathsf W_N(h_r)=2N^{-1}\sum_{i<j}h_r(X_i,X_j)$.  Choose $\delta_0$ so that $8K C_{\mathrm u}\kappa_*t_0\delta_0\le1/2$.  Then
\begin{equation}
 4K C_{\mathrm u}\kappa_*(2t)\frac{m_*}{N}
 \le8K C_{\mathrm u}\kappa_*t_0\delta\le\frac12,
 \qquad
 \E_n\big[e^{t\mathsf W_N(h_r)}\big]\le\left(\frac83\right)^M.
 \label{eq:physical-local-effective-coupling}
\end{equation}
This is uniform over all good occupancy vectors.  Their multinomial weights sum to at most one, so \eqref{eq:physical-good-occupancy-sum} yields
\begin{equation}
 \sup_{N\ge1}\sup_{0\le t\le t_0}
 \E\big[\one_{\mathcal A_{N,J,\delta}}e^{t\mathsf W_N(h_r)}\big]
 \le\left(\frac83\right)^M<\infty.
 \label{eq:physical-centered-local-bound}
\end{equation}
Finally, insert this estimate and \eqref{eq:physical-average-exponential-bound} into \eqref{eq:physical-local-holder}.  The result is at most $(8/3)^{(M+1)/2}$, uniformly in $N\ge1$ and $0\le t\le t_0$.  Both $J_0$ and $\delta_0$ are chosen independently of $N$.
\end{proof}

\begin{proof}[Conclusion of Theorem~\ref{thm:principal}]
The case $\be_0=0$ is immediate.  Otherwise, choose $p,q$ as in \eqref{eq:physical-cutoff-holder}, and take $\delta\le\delta_0$ from Lemma~\ref{lem:physical-local-count-bound} with $t_0=q\be_0$.  Next, choose $J\ge J_0$ large enough that $J^{-\d}<\delta$ and \eqref{eq:physical-bad-gibbs-bound} decays exponentially, and set $r=(2J)^{-1}$.  The singular factor is bounded by Lemma~\ref{lem:physical-local-count-bound}, and the bulk factor by \eqref{eq:physical-continuous-bulk-bound}.  Inserting these bounds into \eqref{eq:physical-cutoff-holder} proves \eqref{eq:principal-order-one}.
\end{proof}

\section{A perturbative estimate for centered arrays}
\label{sec:canonical-arrays}
The variables in a local block have different box-uniform laws.  The following estimate requires independence and centering, but not identical distributions.

\begin{lemma}[Canonical-array exponential estimate]
\label{lem:canonical-array-exponential}
There is a universal $C_{\mathrm u}<\infty$ with the following property.  Let $m\ge1$, let $Y_1,\ldots,Y_m$ be independent, and let $h_{ij}$, $i<j$, be real measurable kernels centered in both variables:
\begin{equation}
 \E[h_{ij}(Y_i,y)]=0,\qquad \E[h_{ij}(x,Y_j)]=0
 \quad\text{for almost every }y,x.
 \label{eq:canonical-array-centering}
\end{equation}
Suppose, for some $\kappa>0$, every integer $p\ge1$, and every $i<j$,
\begin{equation}
 \esssup_y\E[|h_{ij}(Y_i,y)|^p]
 +\esssup_x\E[|h_{ij}(x,Y_j)|^p]\le2\kappa^p p!.
 \label{eq:canonical-array-kernel-moments}
\end{equation}
For $S=\sum_{i<j}h_{ij}(Y_i,Y_j)$, every integer $k\ge1$, and real $\theta$, one has
\begin{equation}
 \E[|S|^{2k}]\le(C_{\mathrm u}\kappa m)^{2k}(2k)!,
 \qquad
 \E[e^{\theta S/m}]
 \le\frac2{1-(C_{\mathrm u}\kappa|\theta|)^2}
 \quad\text{if }C_{\mathrm u}\kappa|\theta|<1.
 \label{eq:canonical-array-exponential}
\end{equation}
Zero kernels are allowed.  We will use the following grouped form.  Partition the $m$ variables into $M$ groups, each containing at most $m_*$ variables.  Suppose each group interacts with at most $K-1$ other groups, where $K\ge1$ is an integer: kernels between all remaining pairs of groups vanish.  Then, for $u\in\R$,
\begin{equation}
 \E[e^{uS/m}]\le\left(\frac83\right)^M
 \quad\text{if}\quad
 4K C_{\mathrm u}\kappa|u|\frac{m_*}{m}\le\frac12.
 \label{eq:canonical-grouped-exponential}
\end{equation}
In both bounds, $m$ is the total number of variables, and the normalization $S/m$ is that of $\mathsf W_N$ in \eqref{eq:homogeneous-partition-function} with $m=N$ particles.  The one-block bound \eqref{eq:canonical-array-exponential} needs $\kappa|\theta|$ small.  The grouped bound instead needs $\kappa|u|\,m_*/m$ small: it trades smallness of the coupling for smallness of the largest group, which is how the occupancy bound $n_a\le\delta N$ enters in Section~\ref{sec:local-singular}.
\end{lemma}

\begin{proof}
We give the counting argument.  Extend the kernels by $h_{ji}(y,x)=h_{ij}(x,y)$.  Independently of the variables $Y_i$, place each index in $I$ independently with probability $1/2$, and put
$D_I=\sum_{i\in I,j\notin I}h_{ij}(Y_i,Y_j)$.
Keep the $Y_i$ fixed and average only over the choice of $I$, denoting this average by $\E_I$.  The transpose convention gives
\begin{equation}
 \begin{aligned}
 D_I&=\sum_{i<j}
 \big(\one_{\{i\in I,\,j\notin I\}}+\one_{\{j\in I,\,i\notin I\}}\big)
 h_{ij}(Y_i,Y_j),\\
 \E_I[D_I]&=\sum_{i<j}\left(\frac14+\frac14\right)h_{ij}(Y_i,Y_j)
 =\frac12 S.
 \end{aligned}
 \label{eq:array-split-identity}
\end{equation}
Each of the two disjoint ways to separate a pair has probability $(1/2)(1/2)=1/4$.  Thus $S=2\E_I[D_I]$ for every fixed sample.  Jensen's inequality, followed by expectation over the $Y_i$, now gives
\begin{equation}
 \E[|S|^{2k}]\le2^{2k}\E_I\E[|D_I|^{2k}].
 \label{eq:array-random-split}
\end{equation}
Fix $I$ and expand $D_I^{2k}$.  Each term is indexed by a row word
$\mathbf i=(i_1,\ldots,i_{2k})$ and a column word
$\mathbf j=(j_1,\ldots,j_{2k})$.
If some label occurs exactly once in either word, the expectation vanishes: condition on all other variables and integrate out that one using \eqref{eq:canonical-array-centering}.  Call a pair of words surviving otherwise; each of its two words then uses at most $k$ distinct labels.

For each label $\ell\in\{1,\ldots,m\}$, let
\begin{equation}
 a_\ell=\#\{v\in\{1,\ldots,2k\}:i_v=\ell\}
 \label{eq:row-multiplicity}
\end{equation}
be the multiplicity of $\ell$ in the row word $\mathbf i$, so that $\sum_\ell a_\ell=2k$.  In a surviving word every $a_\ell\ne1$, so the active labels, those with $a_\ell\ge2$, number at most $k$; labels with $a_\ell=0$ contribute trivial factors to the products below.  Write $T=\prod_{v=1}^{2k}h_{i_vj_v}(Y_{i_v},Y_{j_v})$ and condition on the column variables, $\mathcal C=\sigma(Y_j:j\notin I)$.  Every row index lies in $I$ and every column index outside $I$, so the row variables are independent of $\mathcal C$ and of one another, and the factors sharing the row label $\ell$ form a group of $a_\ell$ functions of $Y_\ell$ alone.  Hence the conditional expectation factors over row labels:
\begin{equation}
 \E[T\mid\mathcal C]
 =\prod_\ell\Phi_\ell(Y_{j_v}:j_v\notin I),
 \qquad
 \Phi_\ell(y)=\E\Big[\prod_{v:\,i_v=\ell}h_{\ell j_v}(Y_\ell,y_{j_v})\Big].
 \label{eq:array-conditional-factor}
\end{equation}
H\"older with all exponents equal to $a_\ell$ and the moment hypothesis \eqref{eq:canonical-array-kernel-moments} give, uniformly in the frozen column values $y$,
\begin{equation}
 |\Phi_\ell(y)|
 \le\prod_{v:\,i_v=\ell}\E\big[|h_{\ell j_v}(Y_\ell,y_{j_v})|^{a_\ell}\big]^{1/a_\ell}
 \le\big(2\kappa^{a_\ell}a_\ell!\big)^{a_\ell/a_\ell}
 =2\kappa^{a_\ell}a_\ell!.
 \label{eq:array-holder-group}
\end{equation}
Since at most $k$ labels are active and $\sum_\ell a_\ell=2k$, the tower property gives $|\E[T]|\le\E\big[|\E[T\mid\mathcal C]|\big]\le\prod_\ell2\kappa^{a_\ell}a_\ell!$, that is,
\begin{equation}
 \left|\E\left[\prod_{v=1}^{2k}
 h_{i_vj_v}(Y_{i_v},Y_{j_v})\right]\right|
 \le 2^k\kappa^{2k}\prod_\ell a_\ell!.
 \label{eq:array-word-moment}
\end{equation}
There are $(2k)!/\prod_\ell a_\ell!$ row words with these multiplicities, so the factorials cancel.  Write $A$ for the number of multiplicity vectors $(a_\ell)_{\ell\le m}$ with $\sum_\ell a_\ell=2k$ and no entry equal to $1$, and $B_k$ for the number of column words of length $2k$ on $m$ labels in which no label occurs exactly once.  Summing absolute values over all surviving pairs of words yields
\begin{equation}
 \E[|D_I|^{2k}]\le2^k\kappa^{2k}(2k)!\,A B_k.
 \label{eq:array-count-reduction}
\end{equation}

It remains to bound $A B_k$.  If $m<2k$, counting all weak compositions gives
$A\le\binom{2k+m-1}{m-1}\le2^{4k}$, and $B_k\le m^{2k}$.
If $m\ge2k$, the active labels of a surviving word form a set of at most $k$ labels.  There are $\sum_{n=0}^k\binom mn\le(em/k)^k$ such sets, and a given set carries at most $4^k$ positive multiplicity vectors and at most $k^{2k}$ words.  Hence
\begin{equation}
 A\le4^k(em/k)^k,\qquad
 B_k\le k^{2k}(em/k)^k,\qquad
 A B_k\le(4e^2)^k m^{2k}.
 \label{eq:array-combinatorial-bound}
\end{equation}
Both cases, together with \eqref{eq:array-random-split}--\eqref{eq:array-count-reduction}, prove the moment bound after enlarging a universal constant.

Finally, $e^x\le2\cosh x$ and the even-moment bound imply
\begin{equation}
 \E[e^{\theta S/m}]
 \le2\sum_{k=0}^\infty
 \frac{|\theta|^{2k}\E[|S|^{2k}]}{m^{2k}(2k)!}
 \le2\sum_{k=0}^\infty(C_{\mathrm u}\kappa|\theta|)^{2k}.
 \label{eq:array-exponential-series}
\end{equation}
Tonelli justifies the nonnegative series.

For the grouped conclusion, draw an edge between each interacting pair of groups.  The graph has degree at most $K-1$.  Greedily color its edges with at most $2K-1$ colors so that edges sharing a group have different colors; place all within-group terms in one further class.  This writes $S=\sum_{c=1}^{2K}\sum_b S_{c,b}$, allowing empty classes.  Within each class, the blocks $S_{c,b}$ use disjoint groups, so they are independent, and at most $M$ of them are nonempty.  A block $S_{c,b}$ involves some number $m_b\le2m_*$ of variables, so \eqref{eq:canonical-array-exponential} applies to it with $m_b$ in place of $m$; taking $\theta=2Kum_b/m$ there gives $\E[e^{2KuS_{c,b}/m}]\le8/3$, because $C_{\mathrm u}\kappa|\theta|=2KC_{\mathrm u}\kappa|u|m_b/m\le4KC_{\mathrm u}\kappa|u|m_*/m\le1/2$.  H\"older over the $2K$ classes then yields
\begin{equation}
 \E[e^{uS/m}]
 \le\prod_{c=1}^{2K}
      \left(\prod_b\E[e^{2KuS_{c,b}/m}]\right)^{1/(2K)}
 \le\left(\frac83\right)^M,
 \label{eq:array-grouped-holder}
\end{equation}
which is \eqref{eq:canonical-grouped-exponential}.
\end{proof}

\appendix
\section{A cell-projection proof of the continuous-kernel bound}
\label{app:continuous-proof}
We prove Lemma~\ref{lem:continuous-entropy-partition}.  The steps are to replace a sample by its histogram and within-cell fluctuation, complete the square in the energy, and estimate the two resulting terms.  The histogram gives a finite-dimensional Gaussian integral; continuity makes the remaining centered kernel small.

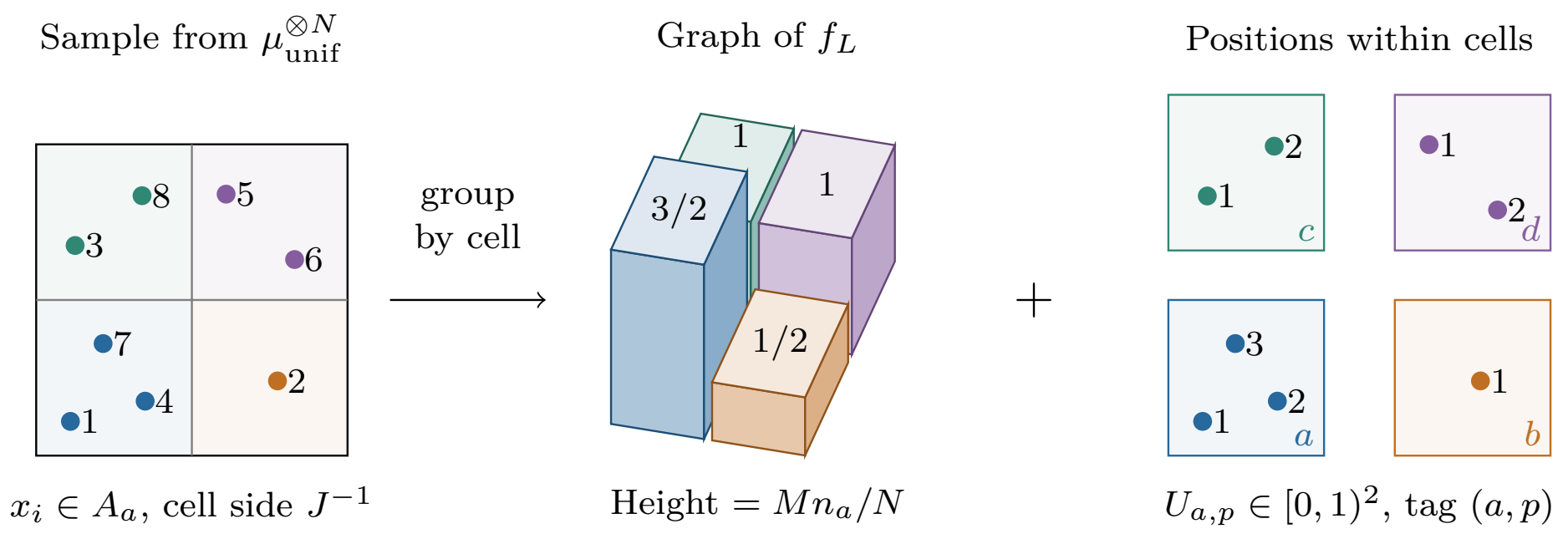
\begin{figure}[htbp]
\centering
\begin{tikzpicture}[x=1.12cm,y=1.12cm,font=\scriptsize,
  particle/.style={circle,inner sep=1.5pt},
  every path/.style={line width=.45pt}]
 \definecolor{cella}{RGB}{39,104,157}
 \definecolor{cellb}{RGB}{191,111,35}
 \definecolor{cellc}{RGB}{48,135,116}
 \definecolor{celld}{RGB}{133,93,158}
 \node at (1.1,2.95) {Sample from $\muunif^{\otimes N}$};
 \node at (5.1,2.95) {Graph of $f_L$};
 \node at (9.35,2.95) {Positions within cells};
 \foreach \x/\y/\col in {0/0/cella,1.1/0/cellb,0/1.1/cellc,1.1/1.1/celld}
   \fill[\col!6] (\x,\y) rectangle +(1.1,1.1);
 \draw (0,0) rectangle (2.2,2.2);
 \draw[gray] (1.1,0)--(1.1,2.2) (0,1.1)--(2.2,1.1);
 \foreach \x/\y/\col/\id in {.242/.242/cella/1,.770/.385/cella/4,.473/.792/cella/7,1.705/.528/cellb/2,.275/1.485/cellc/3,.748/1.837/cellc/8,1.342/1.848/celld/5,1.826/1.386/celld/6} {
   \node[particle,fill=\col] at (\x,\y) {};
   \node[anchor=west,inner sep=2pt] at (\x,\y) {$\id$};
 }
 \node at (1.1,-.35) {$x_i\in A_a$, cell side $J^{-1}$};
 \draw[->] (2.5,1.1)--(3.6,1.1);
 \node[align=center] at (3.05,1.65) {group\\by cell};
 \begin{scope}[shift={(4.02,.20)},x={(.80cm,-.13cm)},
               y={(.37cm,.80cm)},z={(0cm,.92cm)}]
  \foreach \xx/\yy/\hh/\col in {0/1/1/cellc,1/1/1/celld,0/0/1.5/cella,1/0/.5/cellb} {
   \filldraw[fill=\col!38,draw=\col!75!black]
     (\xx+.04,\yy+.04,0)--(\xx+.96,\yy+.04,0)--
     (\xx+.96,\yy+.04,\hh)--(\xx+.04,\yy+.04,\hh)--cycle;
   \filldraw[fill=\col!55,draw=\col!75!black]
     (\xx+.96,\yy+.04,0)--(\xx+.96,\yy+.96,0)--
     (\xx+.96,\yy+.96,\hh)--(\xx+.96,\yy+.04,\hh)--cycle;
   \filldraw[fill=\col!15,draw=\col!75!black]
     (\xx+.04,\yy+.04,\hh)--(\xx+.96,\yy+.04,\hh)--
     (\xx+.96,\yy+.96,\hh)--(\xx+.04,\yy+.96,\hh)--cycle;
  }
  \node at (.5,1.82,1) {$1$};
  \node at (1.5,1.5,1) {$1$};
  \node at (.5,.5,1.5) {$3/2$};
  \node at (1.5,.5,.5) {$1/2$};
 \end{scope}
 \node at (5.1,-.35) {Height $=Mn_a/N$};
 \node[font=\large] at (7.05,1.1) {$+$};
 \foreach \x/\y/\col/\tag in {8/0/cella/a,9.6/0/cellb/b,8/1.45/cellc/c,9.6/1.45/celld/d} {
   \fill[\col!6] (\x,\y) rectangle +(1.1,1.1);
   \draw[\col] (\x,\y) rectangle +(1.1,1.1);
   \node[anchor=south east,inner sep=2pt,text=\col] at (\x+1.1,\y) {$\tag$};
 }
 \foreach \x/\y/\col/\id in {8.242/.242/cella/1,8.770/.385/cella/2,8.473/.792/cella/3,10.205/.528/cellb/1,8.275/1.835/cellc/1,8.748/2.187/cellc/2,9.842/2.198/celld/1,10.326/1.736/celld/2} {
   \node[particle,fill=\col] at (\x,\y) {};
   \node[anchor=west,inner sep=2pt] at (\x,\y) {$\id$};
 }
 \node at (9.35,-.35) {$U_{a,p}\in[0,1)^2$, tag $(a,p)$};
\end{tikzpicture}
\caption{A schematic with $M=4$, $N=8$.  The sample is split into its histogram $f_L$, whose column heights are the cell densities $Mn_a/N$ (height $1$ is uniform), and the positions of the particles within their cells, which carry the signed measure $\nu_N$.  Writing $x_{a,p}=J^{-1}(a+U_{a,p})$ with offsets $U_{a,p}\in[0,1)^\d$, the two right-hand panels recover the sample up to a permutation; each box on the right is a copy of the unit square, with its points numbered within that cell.  Under the uniform product law, conditional on the counts, the offsets are independent and uniform; this is the retagging used in \eqref{eq:partition-occupancy-sum}.}
\label{fig:occupancy-and-offsets}
\end{figure}

\textbf{1. Histogram and energy projection.}
Use the partitions and box averages of Section~\ref{sec:setting}, choosing the resolution $J$ here for the fixed kernel $K$.  In the application $K=b_r$, the cutoff $r$ stays fixed while this auxiliary grid is refined.  With the corresponding labels and counts, set
\begin{equation}
 \mu_N=\frac1N\sum_i\delta_{X_i},\qquad
 f_L=\frac MN\sum_a n_a\one_{A_a},\qquad
 \nu_N=\mu_N-f_L\muunif.
 \label{eq:empirical-histogram-splitting}
\end{equation}
Thus $\nu_N(A_a)=0$.  Figure~\ref{fig:occupancy-and-offsets} shows this decomposition.

For this proof write $T=T_K$ and
$\ip{\sigma,T\tau}=\iint K(x,y)\,\dd\sigma(x)\dd\tau(y)$,
identifying a density with its measure when needed.  Continuity of $K$ makes this form finite on finite signed measures, and positivity extends from densities by approximation.  On the mean-zero cell-constant space $V_M$, define
\begin{equation}
 P_M(\nu_N)=\underset{h\in V_M}{\operatorname{argmin}}\,
                  \ip{h-\nu_N,T(h-\nu_N)}.
 \label{eq:explicit-histogram-square-shift}
\end{equation}
If the form has null directions, take the minimizer modulo them.  Cauchy--Schwarz for the nonnegative form shows that the mixed functional vanishes on these directions.  Projection orthogonality gives
\begin{equation}
 \ip{h,T(\nu_N-P_M(\nu_N))}=0
 \qquad(h\in V_M).
 \label{eq:histogram-square-shift}
\end{equation}
Put $z_N=f_L-1+P_M(\nu_N)$.  Expanding the energy and using \eqref{eq:histogram-square-shift} gives the completed square
\begin{equation}
 \begin{aligned}
 \ip{\mu_N-\muunif,T(\mu_N-\muunif)}
 &=\ip{z_N,Tz_N}+\ip{\nu_N,T\nu_N}
                   -\ip{P_M(\nu_N),TP_M(\nu_N)}\\
 &\le\ip{z_N,Tz_N}+\ip{\nu_N,T\nu_N}.
 \end{aligned}
 \label{eq:direct-square-completion}
\end{equation}

\textbf{2. Entropy bounds the completed square.}
For $h\in V_M$, the entropy variational formula and \eqref{eq:continuous-entropy-assumption} imply
\begin{equation}
 \begin{aligned}
 \log\int e^{Th}\,\dd\muunif
 &=\sup_\rho\left\{\ip{h,T(\rho-1)}-\Hr(\rho\muunif\mid\muunif)\right\}\\
 &\le\sup_\rho\left\{\ip{h,T(\rho-1)}
                         -\frac B2\ip{\rho-1,T(\rho-1)}\right\}
 \le\frac1{2B}\ip{h,Th}.
 \end{aligned}
 \label{eq:dual-log-hls}
\end{equation}
The last step completes the square in $\rho-1-h/B$.  Also, for every $h\in V_M$,
\begin{equation}
 \ip{h,Tz_N}=\int Th\,\dd\mu_N
            =\frac1N\sum_iTh(X_i),
 \label{eq:coarse-gaussian-pairing}
\end{equation}
by the definition of $\nu_N$ and projection orthogonality.

Choose coordinates on the nonnull part of $V_M$, of dimension $q_M\le M-1$, and let $\Sigma$ be the matrix of its energy form.  For a Gaussian coordinate vector $Z$ with covariance $\Sigma^{-1}$, let $\mathsf Z$ denote the corresponding cell-constant function.  The elementary identity
$\E_Z[e^{v^\top\Sigma Z}]=e^{v^\top\Sigma v/2}$ and \eqref{eq:coarse-gaussian-pairing} give, for $0\le t<B$,
\begin{equation}
 \begin{aligned}
 \E\big[e^{\frac{tN}2\ip{z_N,Tz_N}}\big]
 &=\E_Z\left[\left(\int e^{\sqrt{t/N}\,T\mathsf Z}
                           \,\dd\muunif\right)^N\right]\\
 &\le\E_Z\big[e^{\frac{t}{2B}\ip{\mathsf Z,T\mathsf Z}}\big]\\
 &=\frac1{(2\pi)^{q_M/2}\sqrt{\det(\Sigma^{-1})}}
   \int_{\R^{q_M}}e^{-\frac12(1-t/B)z^\top\Sigma z}\,\dd z\\
 &=(1-t/B)^{-q_M/2}.
 \end{aligned}
 \label{eq:coarse-partition-estimate}
\end{equation}
Here, the inequality is \eqref{eq:dual-log-hls} applied to
$h=\sqrt{t/N}\,\mathsf Z$.  Tonelli justifies exchanging the Gaussian and particle integrals.  The bound is independent of $N$ for the fixed grid.

\textbf{3. Continuity makes the remainder small.}
Set
\begin{equation}
 r_J=(I-Q_J)_x(I-Q_J)_yK.
 \label{eq:general-cell-centering}
\end{equation}
This subtracts the box average in each variable and adds back the joint average.  If $\omega_K(\varepsilon)$ is the modulus of continuity when both arguments move by at most $\varepsilon$, then
\begin{equation}
 \|r_J\|_\infty\le4\omega_K(\sqrt\d/J)\longrightarrow0.
 \label{eq:continuous-small-remainder}
\end{equation}
The kernel $r_J$ differs from $K$ by terms that are cell constant in $x$ or in $y$, and $\nu_N$ integrates cell-constant functions to zero; hence $\ip{\nu_N,T\nu_N}=\ip{\nu_N,T_{r_J}\nu_N}$.  Moreover, $r_J$ has zero average in each variable over every box, so $T_{r_J}$ annihilates the cell-constant density $f_L$, and therefore
$\ip{\nu_N,T\nu_N}=\ip{\mu_N,T_{r_J}\mu_N}=N^{-2}\sum_{i,j}r_J(X_i,X_j)$.
Removing equal particle indices in \eqref{eq:direct-square-completion}, we obtain
\begin{equation}
 \mathsf W_N(K)\le N\ip{z_N,Tz_N}
                  +\mathsf W_N(r_J)+5\|K\|_\infty.
 \label{eq:diagonal-correction-bound}
\end{equation}
Indeed, the diagonal correction is
$N^{-1}\sum_i[r_J(X_i,X_i)-K(X_i,X_i)]$, whose absolute value is at most $5\|K\|_\infty$.

Choose conjugate $p,q>1$ with $p\be_0<B$.  Since $r_J$ is centered in each variable under $\muunif$ and bounded by $\|r_J\|_\infty$, Lemma~\ref{lem:canonical-array-exponential} applies with $\kappa=\|r_J\|_\infty$ and $\theta=q\be$; by \eqref{eq:continuous-small-remainder}, a sufficiently fine fixed grid gives
$\sup_{N,\,0\le\be\le\be_0}\E[e^{q\be\mathsf W_N(r_J)/2}]\le8/3$.
H\"older and \eqref{eq:coarse-partition-estimate} now yield
\begin{equation}
 \E\big[e^{\frac\be2\mathsf W_N(K)}\big]
 \le e^{5\be_0\|K\|_\infty/2}
      (1-p\be_0/B)^{-q_M/(2p)}(8/3)^{1/q}.
 \label{eq:continuous-proof-conclusion}
\end{equation}
This proves Lemma~\ref{lem:continuous-entropy-partition}.



\begin{thebibliography}{dlPnMS95}

\bibitem[Bec93]{Beckner1993}
William Beckner, \emph{Sharp {S}obolev inequalities on the sphere and the
  {M}oser--{T}rudinger inequality}, Annals of Mathematics. Second Series
  \textbf{138} (1993), no.~1, 213--242.

\bibitem[BJW23]{BreschJabinWang2023}
Didier Bresch, Pierre-Emmanuel Jabin, and Zhenfu Wang, \emph{Mean field limit
  and quantitative estimates with singular attractive kernels}, Duke
  Mathematical Journal \textbf{172} (2023), no.~13, 2591--2641.

\bibitem[CdCRS25]{ChodrondeCourcelRosenzweigSerfaty2025}
Antonin Chodron~de Courcel, Matthew Rosenzweig, and Sylvia Serfaty, \emph{The
  attractive log gas: {S}tability, uniqueness, and propagation of chaos},
  Communications of the American Mathematical Society \textbf{5} (2025),
  695--773.

\bibitem[CL92]{CarlenLoss1992}
Eric~A. Carlen and Michael Loss, \emph{Competing symmetries, the logarithmic
  {HLS} inequality and {O}nofri's inequality on {$S^n$}}, Geometric and
  Functional Analysis \textbf{2} (1992), no.~1, 90--104.

\bibitem[DGR26a]{DelgadinoGvalaniRosenzweigAttractiveI2026}
Matias~G. Delgadino, Rishabh Gvalani, and Matthew Rosenzweig, \emph{Attractive
  log gas {I}: partition functions and equilibrium fluctuations}, In
  preparation, 2026.

\bibitem[DGR26b]{DelgadinoGvalaniRosenzweig2026}
\bysame, \emph{Sharp mean-field estimates for diffusive log/{R}iesz gases in
  the {H}ilbert--{S}chmidt regime}, 2026, Version 1, 2 September 2026.

\bibitem[dlPMS95]{delaPenaMontgomerySmith1995}
V\'ictor~H. de~la Pe\~na and Stephen~J. Montgomery-Smith, \emph{Decoupling
  inequalities for the tail probabilities of multivariate {$U$}-statistics},
  The Annals of Probability \textbf{23} (1995), no.~2, 806--816.

\bibitem[Gro26]{Grotto2026}
Francesco Grotto, \emph{Fluctuations of point vortex ensembles at small
  negative inverse temperature}, 2026.

\bibitem[Hoe48]{Hoeffding1948}
Wassily Hoeffding, \emph{A class of statistics with asymptotically normal
  distribution}, The Annals of Mathematical Statistics \textbf{19} (1948),
  no.~3, 293--325.

\bibitem[LL07]{LinLucia2007}
Chang-Shou Lin and Marcello Lucia, \emph{One-dimensional symmetry of periodic
  minimizers for a mean field equation}, Annali della Scuola Normale Superiore
  di Pisa. Classe di Scienze. Serie V \textbf{6} (2007), no.~2, 269--290.

\bibitem[LR26]{LeiRosenzweig2026}
Yuhao Lei and Matthew Rosenzweig, \emph{Sharp {B}eckner--{O}nofri inequalities
  on the flat torus: coefficients, defects, and phase transitions}, Work in
  preparation; working title, 2026.

\bibitem[Oki08]{Okikiolu2008}
Kate Okikiolu, \emph{Extremals for logarithmic {H}ardy--{L}ittlewood--{S}obolev
  inequalities on compact manifolds}, Geometric and Functional Analysis
  \textbf{17} (2008), no.~5, 1655--1684.

\bibitem[Ons49]{Onsager1949}
Lars Onsager, \emph{Statistical hydrodynamics}, Il Nuovo Cimento \textbf{6}
  (1949), no.~Suppl. 2, 279--287.

\bibitem[RW25]{RenWang2025}
Zhenjie Ren and Songbo Wang, \emph{Size of chaos for {G}ibbs measures of mean
  field interacting diffusions}, Probability Theory and Related Fields (2025).
  \href{https://doi.org/10.1007/s00440-025-01435-z}{doi:10.1007/s00440-025-01435-z}.
  Accepted version: \href{https://arxiv.org/abs/2411.14236v3}{arXiv:2411.14236v3}.

\bibitem[WZ26a]{WangZhaoAttractive2026}
Zhenfu Wang and Xianliang Zhao, \emph{Full-subcritical partition-function
  asymptotics for the attractive planar {C}oulomb gas at log--{HLS}
  equilibrium}, Work in preparation, 2026.

\bibitem[WZ26b]{WangZhao2026}
\bysame, \emph{Uniform partition-function estimates for {C}oulomb modulated
  energy at all positive temperatures}, 2026, Version 1, 1 September 2026.

\end{thebibliography}
\end{document}